\documentclass[a4paper,12pt]{amsart}
\usepackage{amssymb,amsthm}

\begin{document}

\newcommand{\opp}{\bowtie }
\newcommand{\pos}{\text {\rm pos}}
\newcommand{\supp}{\text {\rm supp}}
\newcommand{\End}{\text {\rm End}}
\newcommand{\diag}{\text {\rm diag}}
\newcommand{\Lie}{\text {\rm Lie}}
\newcommand{\Ad}{\text {\rm Ad}}
\newcommand{\car}{\mathcal R}
\newcommand{\Tr}{\rm Tr}
\newcommand{\Spec}{\text{\rm Spec}}

\def\ge{\geqslant}
\def\le{\leqslant}
\def\a{\alpha}
\def\b{\beta}
\def\c{\chi}
\def\g{\gamma}
\def\G{\Gamma}
\def\d{\delta}
\def\D{\Delta}
\def\L{\Lambda}
\def\e{\epsilon}
\def\et{\eta}
\def\io{\iota}
\def\o{\omega}
\def\p{\pi}
\def\ph{\phi}
\def\ps{\psi}
\def\r{\rho}
\def\s{\sigma}
\def\t{\tau}
\def\th{\theta}
\def\k{\kappa}
\def\l{\lambda}
\def\z{\zeta}
\def\v{\vartheta}
\def\x{\xi}
\def\i{^{-1}}

\def\mapright#1{\smash{\mathop{\longrightarrow}\limits^{#1}}}
\def\mapleft#1{\smash{\mathop{\longleftarrow}\limits^{#1}}}
\def\mapdown#1{\Big\downarrow\rlap{$\vcenter{\hbox{$\scriptstyle#1$}}$}}

\def\ca{\mathcal A}
\def\cb{\mathcal B}
\def\cc{\mathcal C}
\def\cd{\mathcal D}
\def\ce{\mathcal E}
\def\cf{\mathcal F}
\def\cg{\mathcal G}
\def\ch{\mathcal H}
\def\ci{\mathcal I}
\def\cj{\mathcal J}
\def\ck{\mathcal K}
\def\cl{\mathcal L}
\def\cm{\mathcal M}
\def\cn{\mathcal N}
\def\co{\mathcal O}
\def\cp{\mathcal P}
\def\cq{\mathcal Q}
\def\car{\mathcal R}
\def\cs{\mathcal S}
\def\ct{\mathcal T}
\def\cu{\mathcal U}
\def\cv{\mathcal V}
\def\cw{\mathcal W}
\def\cz{\mathcal Z}
\def\cx{\mathcal X}
\def\cy{\mathcal Y}

\def\tz{\tilde Z}
\def\ty{\tilde Y}
\def\tl{\tilde L}
\def\tc{\tilde C}
\def\ta{\tilde A}
\def\tx{\tilde X}

\def \us {\underline{*}}
\def \os {\overline{*}}
\def \trir {\triangleright}
\def \tril {\triangleleft}
\def \Gdia {G_{{\rm diag}}}
\newtheorem{thm}{Theorem}
\newtheorem{lem}{Lemma}
\newtheorem*{rmk}{Remark}
\newtheorem{prop}{Proposition}
\newtheorem{cor}{Corollary}
\newtheorem*{th1}{Lemma 6$'$}

\title[]{A subalgebra of $0$-Hecke algebra}
\author{Xuhua He}
\address{Department of Mathematics, The Hong Kong University of Science and Technology,
Clear Water Bay, Kowloon, Hong Kong}
\email{maxhhe@ust.hk}
\thanks{The author is partially supported by (USA) NSF grant DMS 0700589 and (HK) RGC grant DAG08/09.SC03.}

\begin{abstract}
Let $(W, I)$ be a finite Coxeter group. In the case where $W$ is a Weyl group, Berenstein and Kazhdan in \cite{BK} constructed a monoid structure on the set of all subsets of $I$ using unipotent $\chi$-linear bicrystals. In this paper, we will generalize this result to all types of finite Coxeter groups (including non-crystallographic types). Our approach is more elementary, based on some combinatorics of Coxeter groups. Moreover, we will calculate this monoid structure explicitly for each type.
\end{abstract}
\maketitle

\subsection*{1.1} Let $W$ be a Coxeter group generated by the simply reflections $s_i$ (for $i \in I$). Let $H$ be the Iwahori-Hecke algebra associated to $W$ with parameter $q=0$, i.e., $H$ is a $\mathbb Q$-algebra generated by $T_{s_i}$ for $s_i \in I$ with relations $T_{s_i}^2=-T_{s_i}$ and the braid relations. The algebra $H$ is called $0$-Hecke algebra. It was introduced by Norton in \cite{N}. Representations of $H$ were later studied in the work of Carter \cite{Ca},  Hivert-Novelli-Thibon \cite{HNT} and etc. More recently, Stembridge \cite{St} used the $0$-Hecke algebra to obtain a new proof for the M\"obius function of the Bruhat order of $W$. 

\subsection*{1.2} Set $T'_{s_i}=-T_{s_i}$. For $w \in W$, we define $T'_w=T'_{s_{i_1}} \cdots T'_{s_{i_k}}$, where $w=s_{i_1} \cdots s_{i_k}$ is a reduced expression of $w$. Tit's theorem implies that $T'_w$ is well defined. Moreover, we have a binary operation $\ast: W \times W \to W$ such that $T'_x T'_y=T'_{x \ast y}$ for any $x, y \in W$. It is easy to see that $(W, \ast)$ is a monoid with unit element $1$. 

Now we state our main theorem.

\begin{thm}
Let $W$ be a finite Coxeter group. For any subset $J \subset I$, let $w_0^J$ be the maximal element in the subgroup generated by $s_j$ (for $j \in J$). Then $\{w_0^J w_0^I; J \subset I\}$ is a commutative submonoid of $(W, \ast)$. In other words, there exists a commutative monoid structure $\star_I$ on the set of subsets of $I$, such that $$T'_{w_0^{J_1} w_0^I} T'_{w_0^{J_2} w_0^I}=T'_{w_0^{J_1 \star_I J_2} w_0^I}.$$ 
\end{thm}

\begin{rmk}
In the case where $W$ is a Weyl group, this result was discovered by Berenstein and Kazhdan in \cite[Proposition 2.30]{BK}. Their approach was based on unipotent $\chi$-linear bicrystals. The proof below is more elementary. It is based on some combinatorial properties of Coxeter groups. In the end, we will calculate the operator $\star_I$ explicitly for each type. 
\end{rmk}

\

Below we introduce a binary operation $\trir: W \times W \to W$ and study some properties about the operations $\ast$ and $\trir$. In the case where $W$ is a finite Weyl group, these properties were proved in \cite{HL} using geometry of flag varieties.

\subsection*{1.3} 
We denote by $l$ the length function and $\le$ the Bruhat order on $W$. 

By \cite[Lemma 1.4(1)]{H}, for $x, y \in W$, the subset $\{u y; u \le x\}$ of $W$ contains a unique minimal element which we denote by $x \trir y$. Moreover, $x \trir y=u' y$ for some $u' \le x$ and $l(x \trir y)=l(y)-l(u')$. We also have that if $s_i x>x$, then $(s_i x) \trir y=\min\{s_i (x \trir y), x \trir y\}$.

There is a similar description for $x \ast y$. 

\begin{lem}
Let $x, y \in W$. Then the subset $\{u v; u \le x, v \le y\}$ contains a unique maximal element, which equals $x \ast y$. Moreover, $x \ast y=u' y=x v'$ for some $u' \le x$ and $v' \le y$ and $l(x \ast y)=l(u')+l(y)=l(x)+l(v')$. 
\end{lem}

\begin{rmk}
A slightly weaker version was proved in \cite[Lemma 1.4 (2)]{H}. The proof here is similar to loc.cit.
\end{rmk}

By definition, $s_i \ast w=\max\{w, s_i w\}$. Now for $a, b \in W$ and $i \in I$ with $s_i a>a$, we have that $$T'_{(s_i a) \ast b}=T'_{s_i a} T'_b=T'_{s_i} T'_a T'_b=T'_{s_i} T'_{a \ast b}.$$ Hence $(s_i a) \ast b=\max\{a \ast b, s_i (a \ast b)\}$ if $s_i a>a$. 

We only prove that $\{u v; u \le x, v \le y\}$ contains a unique maximal element which equals $x \ast y$ and $x \ast y=u' y$ for some $u' \le x$ with $l(x \ast y)=l(u')+l(y)$. 

We argue by induction on $l(x)$. For $l(x)=0$, the statement is clear. Assume that $l(x)>0$ and that the statement holds for all $x'$ with $l(x')<l(x)$. Then there exists $i \in I$ such that $s_i x<x$.  By induction hypothesis, the subset $\{u v; u \le s_i x, v \le y\}$ contains a unique maximal element $(s_i x) \ast y$ and there exists $u_1 \le s_i x$ such that $(s_i x) \ast y=u_1 y$ and $l(u_1 y)=l(u_1)+l(y)$. 

Set $z=x \ast y$. Then $z=\max\{u_1 y, s_i u_1 y\}$ and $z y \i=u_1$ or $s_i u_1$. In either case, we have that $z y \i \le x$. If $z y \i=u_1$, then we already know that $l(z)=l(z y \i)+l(y)$. If $z y \i=s_i u_1$, then $s_i u_1 y>u_1 y$ and $$l(s_i u_1 y)=l(u_1 y)+1=l(u_1)+l(y)+1 \ge l(s_i u_1)+l(y) \ge l(s_i u_1 y).$$ Thus we must have that $l(s_i u_1)=l(u_1)+1$ and $l(z)=l(z y \i)+l(y)$. 

Now for any $u \le x$ and $v \le y$. By \cite[Corollary 2.5]{L}, $u \le s_i x$ or $s_i u \le s_i x$. By the definition of $(s_i x) \ast y$, we have that $u v$ or $s_i u v$ is less than or equal to $(s_i x) \ast y=u_1 y \le z$. By \cite[Corollary 2.5]{L}, $u v \le z$. Therefore $z$ is the unique maximal element in the subset $\{u v; u \le x, v \le y\}$.  \qed

\begin{cor}
Let $x' \le x$ and $y' \le y$. Then $x' \ast y' \le x \ast y$.
\end{cor}

\begin{lem}
Let $x' \ge x$ and $y' \le y$, then $x' \trir y' \le x \trir y$.
\end{lem}

By definition, $x' \trir y' \le x \trir y'$. Now we prove that $x \trir y' \le x \trir y$ by induction on $l(x)$. 

For $l(x)=0$, $x \trir y'=y' \le y=x \trir y$. Now assume that $l(x)>0$ and that $x_1 \trir y' \le x_1 \trir y$ for any $x_1$ with $l(x_1)<l(x)$. Then there exists $i \in I$ such that $s_i x<x$.  By induction hypothesis, $(s_i x) \trir y' \le (s_i x) \trir y$. 

By \cite[Corollary 2.5]{L}, $x \trir y'=\min\{(s_i x) \trir y', s_i ((s_i x) \trir y')\} \le (s_i x) \trir y, s_i ((s_i x) \trir y)$. Hence $x \trir y' \le \min\{(s_i x) \trir y, s_i ((s_i x) \trir y)\}=x \trir y$. The statement is proved. \qed

\begin{lem}\label{2}
The action $(W, \ast) \times W \to W$, $(x, y) \mapsto x \trir y$ is a left action of the monoid $(W, \ast)$.
\end{lem}

By definition, $1 \trir x=x$ for any $x \in W$.

Let $x, y, z \in W$. Then there exists $u \le x$ and $v \le y$ such that $y \trir z=v z$ and $x \trir (v z)=u v z$. By definition $u v \le x \ast y$. So $(x \ast y) \trir z \le x \trir (y \trir z)$. 

On the other hand, there exists $w \le x$ such that $x \ast y=w y$ and $l(w y)=l(w)+l(y)$. Then there exists $w' \le w y$ such that $(w y) \trir z=w' z$. Since $l(w y)=l(w)+l(y)$, we may write $w'$ as $w'=w_1 y_1$ for some $w_1 \le w$ and $y_1 \le y$. Thus $(x \ast y) \trir z=w_1 (y_1 z)$. Notice that $w_1 \le x$ and $y_1 \le y$ and $y_1 z \ge y \trir z$. By the previous lemma, $x \trir (y \trir z) \le w_1 (y_1 z)=(x \ast y) \trir z$. The lemma is proved. \qed

\begin{lem} Assume that $W$ is finite. Then $(x \trir y) w_0^I=x \ast (y w_0^I)$.
\end{lem}

By definition, $x \trir y$ is the unique minimal element in $\{u y; u \le x\}$. Hence $(x \trir y) w_0^I$ is the unique maximal element in $\{u y w_0; u \le x\}$. By Lemma 1, this unique maximal element is $x \ast (y w_0^I)$. \qed 

\begin{lem}
Let $J_1, J_2 \subset I$. If $(w_0^{J_1} w_0^I) \ast (w_0^{J_2} w_0^I)=w_0^{J_3} w_0^I$ for some $J_3 \subset I$, then $(w_0^{J_2} w_0^I) \ast (w_0^{J_1} w_0^I)=w_0^{J_3} w_0^I$. In other words,  if $(w_0^{J_1} w_0^I) \trir w_0^{J_2} =w_0^{J_3}$ for some $J_3 \subset I$, then $(w_0^{J_2} w_0^I) \trir w_0^{J_1}=w_0^{J_3}$. 
\end{lem}

This was proved in \cite[Proposition 2.30 (c)]{BK} by applying the anti-automorphism $w \mapsto w_0^I w \i w_0^I$ on $W$. 

\begin{lem}
If $J_2=K \sqcup K'$ with $s_k s_{k'}=s_{k'} s_k$ for any $k \in K$ and $k' \in K'$. Then for any $J_1 \subset I$, we have that $$(w_0^{J_1} w_0^I) \trir w_0^{J_2}=\bigl((w_0^{J_1} w_0^I) \trir w_0^K  \bigr) \bigl((w_0^{J_1} w_0^I) \trir w_0^{K'} \bigr).$$
\end{lem}

We fix a reduced expression $w_0^{J_1} w_0^I=s_{i_1} s_{i_2} \cdots s_{i_n}$ for $i_1, i_2, \cdots, i_n \in I$. Assume that $v \le w_0^{J_1} w_0^I$ with $(w_0^{J_1} w_0^J) \trir w_0^{J_2}=v w_0^{J_2}$. Then $v \le w_0^{J_2}$. Hence $v \in W_{J_2}$. Let $v=s_{i_{l_1}} \cdots s_{i_{l_k}}$ be a reduced subexpression. Then $l_1, \cdots, l_k \in J_2=K \cup K'$. By assumption on $K$ and $K'$, $v=v_1 v_2$ for $v_1 \in W_K$ and $v_2 \in W_{K'}$. So $v w_0^{J_2}=v_1 v_2 w_0^K w_0^{K'}=(v_1 w_0^K) (v_2 w_0^{K'}) \ge \bigl((w_0^{J_1} w_0^I) \trir w_0^K  \bigr) \bigl((w_0^{J_1} w_0^I) \trir w_0^{K'} \bigr)$.

On the other hand, assume that $v_1, v_2 \le w_0^{J_1} w_0^I$ with $(w_0^{J_1} w_0^I) \trir w_0^K=v_1 w_0^K$ and $(w_0^{J_1} w_0^I) \trir w_0^{K'}=v_2 w_0^{K'}$. Then $v_1 \in W_K$ and $v_2 \in W_{K'}$ and there exist $1 \le t_1<t_2<\cdots<t_u \le n$ and $1 \le t'_1<t'_2<\cdots<t'_{u'} \le n$ such that $v_1=s_{i_{t_1}} \cdots s_{i_{t_u}}$ is a reduced subexpression of $v_1$ and $v_2=s_{i_{t'_1}} \cdots s_{i_{t'_{u'}}}$ is a reduced subexpression of $v_2$.

Since $K \cap K'=\emptyset$, $\{t_1, \cdots, t_u\}$ and $\{t'_1, \cdots, t'_{u'}\}$ are disjoint subsets of $\{1, \cdots, n\}$. Let $v$ be the element that corresponds to the subexpression $\{t_1, \cdots, t_u\} \sqcup \{t'_1, \cdots, t'_{u'}\}$. Then it is easy to see that $v=v_1 v_2$. Hence $\bigl((w_0^{J_1} w_0^I) \trir w_0^K  \bigr) \bigl((w_0^{J_1} w_0^I) \trir w_0^{K'} \bigr)=(v_1 w_0^K) (v_2 w_0^{K'})=v_1 v_2 w_0^K w_0^{K'}=v w_0^{J_1} \le (w_0^{J_1} w_0^I) \trir w_0^{J_2}$. \qed

\begin{lem}
Let $J_1, J_2 \subset I$. Then for any $J_1 \subset J'_1$, we have that $$(w_0^{J_1} w_0^I) \trir w_0^{J_2}=(w_0^{J_1} w_0^{J'_1}) \trir \bigl((w_0^{J'_1} w_0^I) \trir w_0^{J_2} \bigr).$$
\end{lem}

Notice that $w_0^{J_1} w_0^I=(w_0^{J_1} w_0^{J'_1}) (w_0^{J_2} w_0^I)=(w_0^{J_1} w_0^{J'_1}) \ast (w_0^{J'_1} w_0^I)$. The lemma follows from Lemma 3. \qed

\

Below is the key lemma.

\begin{lem}\label{4} Assume that $W$ is a irreducible finite Coxeter group and $i, i' \in I$ are end points of the Coxeter graph of $W$. Let $J_1=I-\{i\}$ and $J_2=I-\{i'\}$. Then $(w_0^{J_1} w_0^I) \trir w_0^{J_2}=w_0^{J_3}$ for some $J_3 \subset J_1 \cap J_2$.
\end{lem}

We will prove this lemma in subsection 1.5. The proof is based on a case-by-case checking. We will also use the result to give an explicit description of the operator $\star_I$ for each type. 

Before proving the lemma, we will show that the key lemma implies the main theorem.

\subsection*{1.4 Proof of Theorem 1} By Lemma 5, if $T'_{w_0^{J_1} w_0^I} T'_{w_0^{J_2} w_0^I}=T'_{w_0^{J_3} w_0^I}$, then $T'_{w_0^{J_2} w_0^I} T'_{w_0^{J_1} w_0^I}=T'_{w_0^{J_3} w_0^I}$. Using Lemma 4, we may reformulate the main theorem as follows:

 For any $J_1, J_2 \subset I$, we have that $(w_0^{J_1} w_0^I) \trir w_0^{J_2}=w_0^{J_3}$ for some $J_3 \subset I$. 

We argue by induction on the cardinality of $I$. By Lemma 6 and Lemma 7, it suffices to prove the case where $W$ is irreducible and $J_1, J_2$ are connected in the Coxeter graph. 

It is easy to see that $(w_0^J w_0^I) \trir w_0^I=(w_0^I w_0^I) \trir w_0^J=w_0^J$ and $(w_0^{\emptyset} w_0^I) \trir w_0^J=(w_0^J w_0^I) \trir w_0^{\emptyset}=1=w_0^{\emptyset}$. Now assume that $J_1, J_2$ are proper connected subgraph in the Coxeter graph. Then there exists end points $i, i' \in I$ such that $i \notin J_1$ and $i' \notin J_2$. Set $J'_1=I-\{i\}$ and $J'_2=I-\{i'\}$. Then $J_1 \subset J'_1$ and $J_2 \subset J'_2$. 

By Lemma 7, $(w_0^{J_2} w_0^I) \trir w_0^{J'_1}=(w_0^{J_2} w_0^{J'_2}) \trir \bigl((w_0^{J'_2} w_0^I) \trir w_0^{J'_1} \bigr)$. By Lemma 8, $(w_0^{J'_2} w_0^I) \trir w_0^{J'_1}=w_0^{J_3}$ for some $J_3 \subset J'_1 \cap J'_2$. By induction hypothesis on $W_{J'_2}$, we have that $(w_0^{J_2} w_0^{J'_2}) \trir w_0^{J_3}=w_0^{J_4}$ for some $J_4 \subset J'_1 \cap J_2$. By Lemma 5, $(w_0^{J'_1} w_0^I) \trir w_0^{J_2}=w_0^{J_4}$.

Again by Lemma 7, $(w_0^{J_1} w_0^I) \trir w_0^{J_2}=(w_0^{J_1} w_0^{J'_1}) \trir \bigl((w_0^{J'_1} w_0^I) \trir w_0^{J_2} \bigr)=(w_0^{J_1} w_0^{J'_1}) \trir w_0^{J_4}$. By induction hypothesis on $W_{J'_1}$, we have that $(w_0^{J_1} w_0^{J'_1}) \trir w_0^{J_4}=w_0^{J_5}$ for some $J_5 \subset J_1 \cap J_2$. \qed

\subsection*{1.5 Proof of Lemma \ref{4}} 

We use the same labeling of Coxeter graph as in \cite{Bo}.

For $1 \le a, b \le n$, set $$s_{[a, b]}=\begin{cases} s_a s_{a-1} \cdots s_b, & \text{ if } a \ge b, \\ 1, & \text{ otherwise}. \end{cases}$$

\

{\bf Type $A_n$}

We have that $w_0^{I-\{1\}} w_0^I=s_{[n, 1]} \i$ and $w_0^{I-\{n\}} w_0^I=s_{[n, 1]}$. Hence
\begin{gather*} (w_0^{I-\{1\}} w_0^I) \trir w_0^{I-\{1\}}=s_{[n, 1]} \i \trir w_0^{I-\{1\}}=s_{[n, 2]} \i w_0^{I-\{1\}}=w_0^{I-\{1, 2\}}, \\ (w_0^{I-\{1\}} w_0^I) \trir w_0^{I-\{n\}}=s_{[n, 1]} \i \trir w_0^{I-\{n\}}=s_{[n-1, 1]} \i w_0^{I-\{n\}}=w_0^{I-\{1, n\}}, \\ (w_0^{I-\{n\}} w_0^I) \trir w_0^{I-\{n\}}=s_{[n, 1]} \trir w_0^{I-\{n\}}=s_{[n-1, 1]} w_0^{I-\{n\}}=w_0^{I-\{n-1, n\}}. \end{gather*}

\

{\bf Type $B_n$}

We have that $w_0^{I-\{1\}} w_0^I=s_{[n-1, 1]} \i s_{[n, 1]}$ and $w_0^{I-\{n\}} w_0^I=s_n s_{[n, n-1]} \i \cdots s_{[n, 1]} \i$. Hence
\begin{align*} (w_0^{I-\{1\}} w_0^I) \trir w_0^{I-\{1\}} &=(s_{[n-1, 1]} \i s_{[n, 1]}) \trir w_0^{I-\{1\}}=s_{[n-1, 2]} \i s_{[n, 2]} w_0^{I-\{1\}} \\ &=w_0^{I-\{1, 2\}}, \\ (w_0^{I-\{1\}} w_0^I) \trir w_0^{I-\{n\}} &=(s_{[n-1, 1]} \i s_{[n, 1]}) \trir w_0^{I-\{n\}}=s_{[n-1, 1]} \i \trir (s_{[n, 1]} \trir w_0^{I-\{n\}}) \\ &=s_{[n-1, 1]} \i \trir w_0^{I-\{n-1, n\}}=w_0^{I-\{1, n-1, n\}}, \\ (w_0^{I-\{n\}} w_0^I) \trir w_0^{I-\{n\}} &=(s_n s_{[n, n-1]} \i \cdots s_{[n, 1]} \i) \trir w_0^{I-\{n\}} \\ &=s_{n-1} s_{[n-1, n-2]} \i \cdots s_{[n-1, 1]} \i w_0^{I-\{n\}}=1. \end{align*}

{\bf Type $D_n$}

Set $$\e=\begin{cases} 1, & \text{ if } 2 \nmid n; \\ 0, & \text{ if } 2 \mid n. \end{cases}$$

We have that \begin{gather*} w_0^{I-\{1\}} w_0^I=s_{[n-2, 1]} \i s_{[n, 1]}, \\ w_0^{I-\{n-1\}} w_0^I=s_{n-1} (s_{n-2} s_n) \cdots (s_{[n-2, 2]} \i s_{n-\e}) (s_{[n-2, 1]} \i s_{n-1+\e}), \\ w_0^{I-\{n\}} w_0^I=s_n (s_{n-2} s_{n-1}) \cdots (s_{[n-2, 2]} \i s_{n-1+\e}) (s_{[n-2, 1]} \i s_{n-\e}). \end{gather*}

Hence \begin{align*} (w_0^{I-\{1\}} w_0^I) \trir w_0^{I-\{1\}} &=(s_{[n-2, 1]} \i s_{[n, 1]}) \trir w_0^{I-\{1\}}=s_{[n-2, 2]} \i s_{[n, 2]} w_0^{I-\{1\}} \\ &=w_0^{I-\{1, 2\}}, \\ (w_0^{I-\{1\}} w_0^I) \trir w_0^{I-\{n-1\}} &=(s_{[n-2, 1]} \i s_{[n, 1]}) \trir w_0^{I-\{n-1\}}=s_{[n-2, 1]} \i s_n s_{[n-2, 1]} w_0^{I-\{n-1\}} \\ &=w_0^{I-\{1, n-1, n\}}, \\
(w_0^{I-\{n-1\}} w_0^I) \trir w_0^{I-\{n\}} &=\bigl( s_{n-1} (s_{n-2} s_n) \cdots (s_{[n-2, 2]} \i s_{n-\e}) (s_{[n-2, 1]} \i s_{n-1+\e}) \bigr) \trir w_0^{I-\{n\}} \\ &=s_{n-1} s_{n-2} \cdots s_{[n-2+\e, 2]} \i s_{[n-1-\e, 1]} \i w_0^{I-\{n\}} \\ &=\begin{cases} s_1 s_3 \cdots s_{n-2}, & \text{ if } 2 \nmid n; \\ s_2 s_4 \cdots s_{n-2}, & \text{ if } 2 \mid n. \end{cases} \\ (w_0^{I-\{n\}} w_0^I) \trir w_0^{I-\{n\}} &=\bigl(s_n (s_{n-2} s_{n-1}) \cdots (s_{[n-2, 2]} \i s_{n-1+\e}) (s_{[n-2, 1]} \i s_{n-\e}) \bigr) \trir w_0^{I-\{n\}} \\ &=s_{[n-1, n-2]} \i \cdots s_{[n-1-\e, 2]} \i s_{[n-2+\e, 1]} \i w_0^{I-\{n\}} \\ &=\begin{cases} s_2 s_4 \cdots s_{n-1}, & \text{ if } 2 \nmid n; \\ s_1 s_3 \cdots s_{n-1}, & \text{ if } 2 \mid n. \end{cases} \end{align*}

Applying the automorphism $\s: W \to W$ which exchanges $s_{n-1}$ and $s_n$, we also have that $(w_0^{I-\{1\}} w_0^I) \trir w_0^{I-\{n\}}=w_0^{I-\{1, n-1, n\}}$ and $$(w_0^{I-\{n-1\}} w_0^I) \trir w_0^{I-\{n-1\}}=\begin{cases} (s_2 s_4 \cdots s_{n-3}) s_n, & \text{ if } 2 \nmid n; \\ (s_1 s_3 \cdots s_{n-3}) s_n, & \text{ if } 2 \mid n. \end{cases} .$$

\

For type $E$, set $x=s_4 s_3 s_5 s_4 s_2$.

\

{\bf Type $E_6$}

We have that \begin{gather*} w_0^{I-\{1\}} w_0^I=s_1 s_{[6, 3]} \i x \i s_{[6, 1]} \i, \\ w_0^{I-\{2\}} w_0^I=x \i s_{[6, 1]} \i s_{[5, 1]} x, \\ w_0^{I-\{6\}} w_0^I=s_{[6, 1]} s_{[4, 6]} \i s_{[3, 5]} \i s_2 s_4 s_3 s_1. \end{gather*}

Hence \begin{align*} (w_0^{I-\{1\}} w_0^I) \trir w_0^{I-\{1\}} &=(s_1 s_{[6, 3]} \i x \i s_{[6, 1]} \i) \trir w_0^{I-\{1\}} \\ &=(s_{[6, 3]} \i x \i s_{[6, 2]} \i) w_0^{I-\{1\}}=w_0^{\{2, 4, 5\}},\end{align*} 
\begin{align*} (w_0^{I-\{1\}} w_0^I) \trir w_0^{I-\{2\}} &=(s_1 s_{[6, 3]} \i x \i s_{[6, 1]} \i) \trir w_0^{I-\{2\}} \\ &=(s_1 s_{[6, 3]} \i s_4 s_3 s_5 s_4) \trir \bigl((s_1 s_{[6, 3]} \i) \trir w_0^{I-\{2\}} \bigr) \\ &=(s_1 s_{[6, 3]} \i s_4 s_3 s_5 s_4) \trir w_0^{I-\{1, 2\}} \\ &=s_{[6, 3]} \i s_4 s_3 s_5 s_4 w_0^{I-\{1, 2\}}=s_4 s_6,\end{align*} 
\begin{align*} (w_0^{I-\{1\}} w_0^I) \trir w_0^{I-\{6\}} &=(s_1 s_{[6, 3]} \i x \i s_{[6, 1]} \i) \trir w_0^{I-\{6\}} \\ &=s_1 s_{[5, 3]} \i x \i s_{[5, 1]} \i w_0^{I-\{6\}}=w_0^{\{3, 4, 5\}},\end{align*} 
\begin{align*} (w_0^{I-\{2\}} w_0^I) \trir w_0^{I-\{2\}} &=(x \i s_{[6, 1]} \i s_{[5, 1]} x) \trir w_0^{I-\{2\}} \\ &=(x \i) \trir \bigl((s_1 s_{[6,3]} \i s_{[5,3]} s_1 s_4 s_3 s_5 s_4) \trir w_0^{I-\{2\}} \bigr) \\ &=x \i \trir (s_3 s_5)=1.
\end{align*}

Applying the nontrivial diagram automorphism, we also have that $(w_0^{I-\{6\}} w_0^I) \trir w_0^{I-\{6\}}=w_0^{\{2, 3, 4\}}$ and $(w_0^{I-\{6\}} w_0^I) \trir w_0^{I-\{2\}}=s_1 s_4$.

\

{\bf Type $E_7$}

We have that \begin{gather*} w_0^{I-\{1\}} w_0^I=s_1 s_3 s_4 s_2 s_{[5, 3]} s_1 s_{[6, 2]} s_{[6, 4]} \i s_{[7, 1]} x s_{[6, 3]} s_1, \\ w_0^{I-\{2\}} w_0^I=s_2 s_4 s_3 s_1 s_{[5, 2]} s_{[6, 4]} \i s_{[5, 1]} x s_{[7, 1]} x s_{[6, 3]} s_{[7, 4]} s_2, \\ w_0^{I-\{7\}} w_0^I=s_{[7, 1]} x s_{[6, 3]} s_1 s_{[7, 2]} s_{[7, 4]} \i.\end{gather*}

Hence \begin{align*} (w_0^{I-\{1\}} w_0^I) \trir w_0^{I-\{1\}} &=(s_1 s_3 s_4 s_2 s_{[5, 3]} s_1 s_{[6, 2]} s_{[6, 4]} \i s_{[7, 1]} x s_{[6, 3]} s_1) \trir w_0^{I-\{1\}} \\ &=(s_3 s_4 s_2 s_{[5, 3]}) \trir \bigl((s_{[6, 2]} s_{[6, 4]} \i s_{[7, 2]} x s_{[6, 3]}) \trir w_0^{I-\{1\}} \bigr) \\ &=(s_3 s_4 s_2 s_{[5, 3]}) \trir w_0^{\{2,4,5,7\}}=s_2 s_5 s_7,\end{align*} 
\begin{align*} (w_0^{I-\{1\}} w_0^I) \trir w_0^{I-\{2\}} &=(s_1 s_3 s_4 s_2 s_{[5, 3]} s_1 s_{[6, 2]} s_{[6, 4]} \i s_{[7, 1]} x s_{[6, 3]} s_1) \trir w_0^{I-\{2\}} \\ &=(s_1 s_3 s_4 s_2 s_{[5, 3]} s_1 s_{[6, 2]} s_{[6, 4]} \i) \trir \bigl( (s_{[7, 1]} x s_{[6, 3]} s_1) \trir w_0^{I-\{2\}} \bigr) \\ &=(s_1 s_3 s_4 s_2 s_{[5, 3]} s_1 s_{[6, 2]} s_{[6, 4]} \i) \trir (s_1 s_4 s_{[6, 3]} s_1)=1,\end{align*} 
\begin{align*} (w_0^{I-\{1\}} w_0^I) \trir w_0^{I-\{7\}} &=(s_1 s_3 s_4 s_2 s_{[5, 3]} s_1 s_{[6, 2]} s_{[6, 4]} \i s_{[7, 1]} x s_{[6, 3]} s_1) \trir w_0^{I-\{7\}} \\ &=(s_1 s_3 s_4 s_2 s_{[5, 3]} s_1 s_{[6, 2]} s_{[6, 4]} \i) \trir \bigl((s_{[6, 1]} x s_{[6, 3]} s_1) \trir w_0^{I-\{7\}} \bigr) \\ &=(s_1 s_3 s_4 s_2 s_{[5, 3]} s_1 s_{[6, 2]} s_{[6, 4]} \i) \trir w_0^{I-\{6, 7\}} \\ &=s_1 s_3 s_4 s_2 s_{[5, 3]} s_1 s_{[5, 2]} s_{[5, 4]} \i w_0^{I-\{6, 7\}}=w_0^{\{3,4,5\}}, \end{align*}
\begin{align*} (w_0^{I-\{2\}} w_0^I) \trir w_0^{I-\{2\}} &=(s_2 s_4 s_3 s_1 s_{[5, 2]} s_{[6, 4]} \i s_{[5, 1]} x s_{[7, 1]} x s_{[6, 3]} s_{[7, 4]} s_2) \trir w_0^{I-\{2\}} \\ &=(s_2 s_4 s_3 s_1 s_{[5, 2]} s_{[6, 4]} \i s_{[5, 1]} x) \trir \bigl((s_{[7, 1]} x s_{[6, 3]} s_{[7, 4]}) \trir w_0^{I-\{2\}} \bigr) \\ &=(s_2 s_4 s_3 s_1 s_{[5, 2]} s_{[6, 4]} \i s_{[5, 1]} x) \trir (s_1 s_4 s_6)=1, \end{align*}
\begin{align*} (w_0^{I-\{2\}} w_0^I) \trir w_0^{I-\{7\}} &=(s_2 s_4 s_3 s_1 s_{[5, 2]} s_{[6, 4]} \i s_{[5, 1]} x s_{[7, 1]} x s_{[6, 3]} s_{[7, 4]} s_2) \trir w_0^{I-\{7\}} \\ &=(s_2 s_4 s_3 s_1 s_{[5, 2]} s_{[6, 4]} \i) \trir \bigl( (s_{[5, 1]} x s_{[6, 1]} x s_{[6, 3]} s_{[6, 4]} s_2) \trir w_0^{I-\{7\}} \bigr) \\ &=(s_2 s_4 s_3 s_1 s_{[5, 2]} s_{[6, 4]} \i) \trir w_0^{\{1,3,4,6\}}=1, \end{align*}
\begin{align*} (w_0^{I-\{7\}} w_0^I) \trir w_0^{I-\{7\}} &=(s_{[7, 1]} x s_{[6, 3]} s_1 s_{[7, 2]} s_{[7, 4]} \i) \trir w_0^{I-\{7\}} \\ &=(s_{[6, 1]} x s_{[6, 3]} s_1 s_{[6, 2]} s_{[6, 4]} \i) w_0^{I-\{7\}}=w_0^{\{2,3,4,5\}}.
\end{align*}

\

{\bf Type $E_8$}

We have that \begin{align*} w_0^{I-\{1\}} w_0^I &=s_1 s_3 s_4 s_2 s_{[5,3]} s_1 s_{[6,2]} s_{[6,4]} \i s_{[7,1]} x s_{[6,3]} s_1 s_{[8,1]} x s_{[6,3]} s_1 s_{[7,2]} s_{[6,4]} \i s_{[8,1]} x s_{[6,3]} s_1, \\ w_0^{I-\{2\}} w_0^I &=s_2 s_4 s_3 s_1 s_{[5,2]} s_{[6,4]} \i s_{[5,1]} x s_{[7,1]} x s_{[6,3]} s_{[7,4]} s_2 s_{[8,1]} x s_{[6,3]} s_1 s_{[7,2]} s_{[7,4]} \i s_{[8,1]} x s_{[6,3]} s_{[7,4]} s_2, \\ w_0^{I-\{8\}} w_0^I &=s_{[8,1]} x s_{[6,3]} s_1 s_{[7,2]} s_{[7,4]} \i s_{[8,1]} x s_{[6,3]} s_1 s_{[7,2]} s_{[8,4]} \i.
\end{align*}

Hence \begin{align*} (w_0^{I-\{1\}} w_0^I) \trir w_0^{I-\{1\}} &=(s_1 s_3 s_4 s_2 s_{[5,3]} s_1 s_{[6,2]} s_{[6,4]} \i s_{[7,1]} x s_{[6,3]} s_1) \trir \\ & \quad \bigl( (s_{[8,1]} x s_{[6,3]} s_1 s_{[7,2]} s_{[6,4]} \i s_{[8,1]} x s_{[6,3]} s_1) \trir w_0^{I-\{1\}} \bigr)\\ &=(s_1 s_3 s_4 s_2 s_{[5,3]} s_1 s_{[6,2]} s_{[6,4]} \i s_{[7,1]} x s_{[6,3]} s_1) \trir (s_3 s_5 s_7) \\ &=1, \end{align*}
\begin{align*} (w_0^{I-\{1\}} w_0^I) \trir w_0^{I-\{2\}} &=(s_1 s_3 s_4 s_2 s_{[5,3]} s_1 s_{[6,2]} s_{[6,4]} \i s_{[7,1]} x s_{[6,3]} s_1 s_{[8,1]} x s_{[6,3]} s_1) \trir \\ & \quad \bigl( (s_{[7,2]} s_{[6,4]} \i s_{[8,1]} x s_{[6,3]} s_1) \trir w_0^{I-\{2\}} \bigr) \\ &=(s_1 s_3 s_4 s_2 s_{[5,3]} s_1 s_{[6,2]} s_{[6,4]} \i s_{[7,1]} x s_{[6,3]} s_1 s_{[8,1]} x s_{[6,3]} s_1) \trir (s_1 s_3 s_{[6,3]} s_1) \\ &=1, \end{align*}
\begin{align*} (w_0^{I-\{1\}} w_0^I) \trir w_0^{I-\{8\}} &=(s_1 s_3 s_4 s_2 s_{[5,3]} s_1 s_{[6,2]} s_{[6,4]} \i s_{[7,1]} x s_{[6,3]} s_1) \trir \\ & \quad \bigl((s_{[8,1]} x s_{[6,3]} s_1 s_{[7,2]} s_{[6,4]} \i s_{[8,1]} x s_{[6,3]} s_1) \trir w_0^{I-\{8\}} \bigr) \\ &=(s_1 s_3 s_4 s_2 s_{[5,3]} s_1 s_{[6,2]} s_{[6,4]} \i s_{[7,1]} x s_{[6,3]} s_1) \trir w_0^{\{2,3,4,5,6\}} \\ &=(s_1 s_3 s_4 s_2 s_{[5,3]} s_1 s_{[6,2]} s_{[6,4]} \i) \trir \bigl( (s_{[7,1]} x s_{[6,3]} s_1) \trir w_0^{\{2,3,4,5,6\}} \bigr) \\ &=(s_1 s_3 s_4 s_2 s_{[5,3]} s_1 s_{[6,2]} s_{[6,4]} \i) \trir w_0^{\{3,4,5\}}=1, \end{align*}
\begin{align*} (w_0^{I-\{2\}} w_0^I) \trir w_0^{I-\{2\}} &=(s_2 s_4 s_3 s_1 s_{[5,2]} s_{[6,4]} \i s_{[5,1]} x s_{[7,1]} x s_{[6,3]} s_{[7,4]} s_2 s_{[8,1]} x s_{[6,3]} s_1 s_{[7,2]} s_{[7,4]} \i) \trir \\ & \quad \bigl( (s_{[8,1]} x s_{[6,3]} s_{[7,4]} s_2) \trir w_0^{I-\{2\}}\bigr) \\ &=(s_2 s_4 s_3 s_1 s_{[5,2]} s_{[6,4]} \i s_{[5,1]} x s_{[7,1]} x s_{[6,3]} s_{[7,4]} s_2 s_{[8,1]} x s_{[6,3]} s_1 s_{[7,2]} s_{[7,4]} \i) \trir \\& \quad (s_1 s_4 s_6 s_{[7,3]} s_1) \\ &=1, \end{align*}
\begin{align*} (w_0^{I-\{2\}} w_0^I) \trir w_0^{I-\{8\}} &=(s_2 s_4 s_3 s_1 s_{[5,2]} s_{[6,4]} \i s_{[5,1]} x s_{[7,1]}) \trir \\ & \quad (x s_{[6,3]} s_{[6,4]} s_2 s_{[7,1]} x s_{[6,3]} s_1 s_{[6,2]} s_{[6,4]} \i s_{[7,1]} x s_{[6,3]} s_{[7,4]} s_2 w_0^{I-\{8\}}) \\ &=(s_2 s_4 s_3 s_1 s_{[5,2]} s_{[6,4]} \i s_{[5,1]} x s_{[7,1]}) \trir w_0^{\{1, 5,6\}}=1, \end{align*}
\begin{align*} (w_0^{I-\{8\}} w_0^I) \trir w_0^{I-\{8\}} &=s_{[6,1]} x s_{[6,3]} s_1 s_{[6,2]} s_{[6,4]} \i s_{[7,1]} x s_{[6,3]} s_1 s_{[7,2]} s_{[7,4]} \i w_0^{I-\{8\}} \\&=w_0^{\{2,3,4,5\}}.
\end{align*}

\

{\bf Type $F_4$}

We have that $w_0^{I-\{1\}} w_0^I=s_{[4, 1]} \i s_2 s_3 s_2 s_1 s_2 s_3 s_2 s_{[4, 1]}$ and $w_0^{I-\{4\}} w_0^I=s_{[4, 1]} s_3 s_2 s_3 s_4 s_3 s_2 s_3 s_{[4, 1]} \i$. Hence \begin{gather*} (w_0^{I-\{1\}} w_0^I) \trir w_0^{I-\{1\}}=(s_{[4, 1]} \i s_2 s_3 s_2 s_1 s_2 s_3 s_2 s_{[4, 1]}) \trir w_0^{I-\{1\}}=1, \\ (w_0^{I-\{1\}} w_0^I) \trir w_0^{I-\{4\}}=(s_{[4, 1]} \i s_2 s_3 s_2 s_1 s_2 s_3 s_2 s_{[4, 1]}) \trir w_0^{I-\{4\}}=1.
\end{gather*}

Applying the nontrivial diagram automorphism, we also have that $(w_0^{I-\{4\}} w_0^I) \trir w_0^{I-\{4\}}=1$.

\

{\bf Type $H_3$}

We have that $w_0^{I-\{1\}} w_0^I=s_1 s_2 s_1 s_2 s_3 s_2 s_1 s_2 s_1 s_3 s_2 s_1$ and $w_0^{I-\{3\}} w_0^I=s_3 s_2 s_1 s_2 s_1 s_3 s_2 s_1 s_2 s_3$. Hence it is easy to see that $(w_0^{I-\{1\}} w_0^I) \trir w_0^{I-\{1\}}=(w_0^{I-\{1\}} w_0^I) \trir w_0^{I-\{3\}}=(w_0^{I-\{3\}} w_0^I) \trir w_0^{I-\{3\}}=1$.

\

{\bf Type $H_4$}

We have that $w_0^{I-\{1\}} w_0^I=s_1 s_2 s_1 s_2 s_3 s_2 s_1 s_2 s_1 s_3 s_2 s_1 (s_{[4,1]} s_2 s_1 s_{[3,1]} s_2 s_3)^4 s_3 s_2$ and $w_0^{I-\{4\}} w_0^I=(s_{[4,1]} s_2 s_1 s_{[3,1]} s_2 s_3)^4 s_4$. Hence it is easy to see that $(w_0^{I-\{1\}} w_0^I) \trir w_0^{I-\{1\}}=(w_0^{I-\{1\}} w_0^I) \trir w_0^{I-\{4\}}=(w_0^{I-\{4\}} w_0^I) \trir w_0^{I-\{4\}}=1$.

\

{\bf Type $I_m$}

It is easy to see that $(w_0^{I-\{1\}} w_0^I) \trir w_0^{I-\{1\}}=(w_0^{I-\{1\}} w_0^I) \trir w_0^{I-\{2\}}=(w_0^{I-\{2\}} w_0^I) \trir w_0^{I-\{2\}}=1$.

\subsection*{1.6} Now we will calculate the operator $\star_I$ explicitly for each type. This is based on the previous subsection, the equalities $J_1 \star_I J_2=J_2 \star_I J_1$ (see Lemma 5), $J'_1 \star_I J_2=J'_1 \star_{J_1} (J_1 \star_I J_2)$ for any $J'_1 \subset J_1$ (see Lemma 7) and the inequality $J_1 \star_I J_2 \subset J'_1 \star_I J'_2$ for $J_1 \subset J'_1$ and $J_2 \subset J'_2$ (see Lemma 2). We just list below the cases where $J_1$ and $J_2$ are proper connected subgraph of the Coxeter graph of $W$.  

For the case where $J_1$ and $J_2$ are not necessarily connected and $J'_1, \cdots, J'_l$ (resp. $J''_1, \cdots, J''_{l'}$) are the connected components of $J_1$ (resp. $J_2$), we have that $J_1 \star_I J_2=\sqcup_{1 \le i \le l, 1 \le i' \le l'} (J'_i \star_I J''_{i})$ (see Lemma 6). 

\

{\bf Type $A_n$:} Let $J_1=\{a, a+1, \cdots, n-b\}, J_2=\{a', a'+1, \cdots, n-b'\}$. Then $$J_1 \star_I J_2=\{a+a'-1, a+a', \cdots, n-b-b'\}.$$ 

\

{\bf Type $B_n$:} Let $J_1=\{a, a+1, \cdots, n-b\}, J_2=\{a', a'+1, \cdots, n-b'\}$. Then \[J_1 \star_I J_2=\begin{cases} \{a+a'-1, a+a', \cdots, n\}, & \text{ if } b=b'=0; \\ \{a+a'-1, a+a', \cdots, n-b'-1\}, & \text{ if } b=0, b' \ge 1; \\ \emptyset, & \text{ if } b, b' \ge 1.\end{cases}\]

\

{\bf Type $D_n$:} If $I-J_1=n-1$ or $n$ and $I-J_2=n-1$ or $n$, then $J_1 \star_I J_2$ was already calculated in the previous subsection.

If $\{n-1, n\} \nsubseteq J_1$ and $\{n-1, n\} \nsubseteq J_2$ and $J_1$ or $J_2$ contains at most $n-2$ elements, then $J_1 \star_I J_2=\emptyset$.  

Otherwise, we may assume that $\{n-1, n\} \subset J_1$, i.e. $J_1=\{a, a+1, \cdots, n\}$ for some $a$. If $\{n-1, n\} \subset J_2$, i.e. $J_2=\{a', a'+1, \cdots, n\}$ for some $a'$, then $J_1 \star_I J_2=\{a+a'-1, a+a', \cdots n\}$. 

If $\{n-1, n\} \nsubseteq J_2$, we may assume without loss of generality that $n \notin J_2$, i.e. $J_2=\{a', a'+1, \cdots, n-b'\}$ for some $a', b'$ with $b' \ge 1$. Then $J_1 \star_I J_2=\{a+a'-1, a+a', \cdots, n-b'-1\}$. 

\

{\bf Type $E_6$:} If $2 \notin J_1$ and $2 \notin J_2$, then $J_1 \star_I J_2=\emptyset$.

If $2 \in J_1$ and $2 \notin J_2$, then \[J_1 \star_I J_2=\begin{cases} \{4,6\}, & \text{ if } J_1=I-\{1\}, J_2=I-\{2\}; \\ \{1, 4\}, & \text{ if } J_1=I-\{6\}, J_2=I-\{2\}; \\ \emptyset, & \text{ otherwise}. \end{cases}\]

If $2 \in J_1 \cap J_2$, then \[J_1 \star_I J_2=\begin{cases} \{2, 4, 5\}, & \text{ if } J_1=J_2=I-\{1\}, \\ \{2, 3, 4\}, & \text{ if } J_1=J_2=I-\{6\}, \\ \{3, 4, 5\}, & \text{ if } \{J_1, J_2\}=\{I-\{1\}, J_2=I-\{6\}\}, \\ \{4\}, & \text{ if } \{J_1, J_2\}=\{I-\{1\}, I-\{1, 6\}\}, \\ \{4\}, & \text{ if } \{J_1, J_2\}=\{I-\{6\}, I-\{1, 6\}\}, \\ \emptyset, & \text{ otherwise}. \end{cases}\]

\

{\bf Type $E_7$:} For proper subsets $J_1$ and $J_2$, we have that \[J_1 \star_I J_2=\begin{cases} \{2, 5, 7\}, & \text{ if } J_1=J_2=I-\{1\}, \\ \{2,3,4,5\}, & \text{ if } J_1=J_2=I-\{7\}, \\ \{3,4,5\}, & \text{ if } \{J_1, J_2\}=\{I-\{1\}, I-\{7\}\}, \\ \{4\}, & \text{ if } \{J_1, J_2\}=\{I-\{7\}, I-\{1, 7\}\}, \\ \{4\}, & \text{ if } \{J_1, J_2\}=\{I-\{7\}, I-\{6, 7\}\}, \\ \emptyset, & \text{ otherwise}. \end{cases}\]

\

{\bf Type $E_8$:} For proper subsets $J_1$ and $J_2$, we have that \[J_1 \star_I J_2=\begin{cases} \{2, 3, 4, 5\}, & \text{ if } J_1=J_2=I-\{8\}, \\ \emptyset, & \text{ otherwise}. \end{cases}\]

\

{\bf Type $F$, $H$ and $I$:} For proper subsets $J_1, J_2$ of $I$, we always have that $J_1 \star_I J_2=\emptyset$. 

\subsection*{Acknowledgements}
During my visit in University of Oregon, Berenstein told me the result \cite[Proposition 2.30]{BK} that was discovered in his joint work with Kazhdan and asked me if I can find a combinatorial explanation. It is my pleasure to thank him. Furthermore, I thank the referee for his/her valuable comments and suggestions, especially for the suggestion of using the language of $0$-Hecke algebras and for pointing out some references on this topics. The computations for exceptional groups were done with the aid of CHEVIE package \cite{CH} of \cite{GA}.

\end{document}